\documentclass[aos,MSNbibl,dvips]{arximspdf}

\doi{10.1214/12-AOS988} 
\volume{40}
\issue{2}
\pubyear{2012}
\firstpage{1283}
\lastpage{1284}

\makeatletter
\newtheorem{theorem}{Theorem}
\newcommand{\E}{\mathrm{E}}
\newcommand{\Var}{\operatorname{Var}}
\newcommand{\trace}{\operatorname{tr}}
\makeatother

\begin{document}
\begin{frontmatter}

\title{Correction on\\
Moments of minors of Wishart matrices}\vspace*{6pt}
\runtitle{Correction}
\pdftitle{Correction on
Moments of minors of Wishart matrices by M. Drton and A. Goia}

\textit{Ann. Statist.} \textbf{36} (2008)
2261--2283

\begin{aug}
\author[A]{\fnms{Mathias} \snm{Drton}\corref{}\ead[label=e1]{drton@galton.uchicago.edu}}
\and
\author[B]{\fnms{Aldo} \snm{Goia}\ead[label=e2]{aldo.goia@eco.unipmn.it}}
\runauthor{M. Drton and A. Goia}
\affiliation{University of Chicago and University of East Piedmont}
\address[A]{Department of Statistics \\
University of Chicago\\
Chicago, Illinois  60637\\
USA\\
\printead{e1}} 
\address[B]{Dipartimento di Studi per l'Economia \\
\quad e l'Impresa \\
University of East Piedmont \\
Novara\\
Italy\\
\printead{e2}}
\end{aug}

\received{\smonth{3} \syear{2012}}



\end{frontmatter}

\def\thetheorem{5.7}

Theorem 5.7 in \cite{drton} gives a formula for the variance of a minor
(i.e., a subdeterminant) of a Wishart random matrix. The formula has to
be corrected as follows.

\begin{theorem}
\label{thm:var-gen}  Let $I,J\in{r\brace m}$ have intersection
$C:=I\cap J$ of  cardinality $c=|C|=|I\cap J|$. Define $\bar
I=I\setminus(I\cap J)$,  $\bar J=J\setminus(I\cap J)$ and $\bar I\bar
J=\bar I \cup\bar J$.  Then the minor
$\det(S_{I\times J})=\det(S_{IJ})$ of the Wishart matrix $S\sim\mathcal{W}%
_{r}(n,\Sigma)$ has  variance
\begin{eqnarray*}
&&\Var[\det(  S_{IJ})]
\\
&&\qquad= \det(\Sigma_{IJ})^{2} \frac{n!}{(n-m)!} \biggl[
\frac{(n+2)!}{(n+2-m)!} - \frac{n!}{(n-m)!} \biggr]  \\
&&\quad\qquad{} + \det(\Sigma_{C\times C})^{2} \det(\bar\Sigma_{\bar I\bar
J\times\bar I\bar
J}) \frac{(n+2)!}{(n+2-c)!} \cdot\frac{n!}{(n-m)!} \\
&&\qquad\qquad{}\times \Biggl[  \sum_{k=0}^{m-c-1}
(m-c-k)!\frac{(n+2-c)!}{(n+2-c-k)!}(-1)^{k}
\trace\bigl\{  (  \bar{\Sigma}_{\bar{I}\bar{J}}\bar{\Sigma}^{\bar{I}%
\bar{J}}) ^{(k)} \bigr\}  \Biggr] ,
\end{eqnarray*}
where $ \bar\Sigma=\Sigma_{([r]\setminus C)\times([r]\setminus C)}-
\Sigma_{([r]\setminus C)\times C}\Sigma_{C\times C}^{-1}\Sigma_{C\times
([r]\setminus C)}$.
\end{theorem}

\begin{pf}
Define $\bar S$ in analogy to $\bar \Sigma$.  Since $\det(S_{I J}) =
\det(S_{C\times C})\det(\bar S_{\bar I\times \bar J})$, and $S_{C\times
C}$ and $\bar S_{\bar I\times \bar J}$ are independent (Lemma~5.2),
\begin{eqnarray*}
\Var[  \det(  S_{IJ})  ]   &=&\bigl(  \Var[ \det( S_{CC}) ] +\E[ \det( S_{CC})
] ^{2}\bigr) \Var[ \det( \bar{S}_{\bar
{I}\bar{J}})  ] \\
&&{}+\Var[  \det(  S_{CC})  ]  \E[ \det( \bar{S}_{\bar{I}\bar{J}}) ] ^{2}.
\end{eqnarray*}
The claim now follows from Corollary~4.2 and Propositions~5.1 and
5.5.
\end{pf}

\def\bibname{Reference}


\printaddresses

\end{document}